\documentstyle[10pt]{article}
\title{Approximate Unitary Equivalence of Homomorphisms
          from ${\cal O}_{\infty}$
   \thanks{Research partially supported by NSF grants DMS 93-01082
  (H. Lin) and DMS 91-06285 and DMS 94-00904 (N. C. Phillips),
  and by the Fields Institute (both authors). \protect\\
   AMS 1991 subject classification numbers: Primary 46L35;
   Secondary 19K99, 46L80, 46M40. \protect\\
 }  }
\author{Huaxin Lin \\ and \\ N. Christopher Phillips \\ \\
  Department of Mathematics \\ University
    of Oregon \\ Eugene OR 97403-1222
		\vspace{1.5in}}
\date{}
\begin{document}
\maketitle
\newcommand{\beq}{\begin{equation}}
\newcommand{\eeq}{\end{equation}}
\newcommand{\bd}{\begin{displaymath}}
\newcommand{\ed}{\end{displaymath}}
\newcommand{\ben}{\begin{enumerate}}
\newcommand{\een}{\end{enumerate}}
\newcommand{\bde}{\begin{description}} 
\newcommand{\ede}{\end{description}}
\newcommand{\beqr}{\begin{eqnarray*}}
\newcommand{\eeqr}{\end{eqnarray*}}
\newcommand{\hs}[1]{\hspace{#1}}
\newcommand{\vs}[1]{\vspace{#1}}
\newcommand{\bc}{\begin{center}}
\newcommand{\ec}{\end{center}}
\newcommand{\bv}{\begin{verbatim}}
\newcommand{\ev}{\end{verbatim}}

\newcommand{\af}{\alpha}

\newcommand{\half}{\frac{1}{2}}
\newcommand{\p}{\partial}
\newcommand{\bt}{\beta}
\newcommand{\gm}{\gamma}
\newcommand{\dt}{\delta}
\newcommand{\ep}{\epsilon}
\newcommand{\zt}{\zeta}
\newcommand{\et}{\eta}
\newcommand{\tta}{\theta}
\newcommand{\ld}{\lambda}
\newcommand{\sm}{\sigma}
\newcommand{\vph}{\varphi}
\newcommand{\kp}{\kappa}
\renewcommand{\i}{\subset}

\newcommand{\s}[1]{\noindent {\bf EXERCISE~#1}{}}

\setcounter{section}{-1}
\newtheorem{thm}{Theorem}[section]
\newtheorem{lem}[thm]{Lemma}
\newtheorem{prop}[thm]{Proposition}
\newtheorem{dfn}[thm]{Definition}
\newtheorem{cor}[thm]{Corollary}
\newtheorem{conj}[thm]{Conjecture}
\newtheorem{conv}[thm]{Convention}
\newtheorem{rmk}[thm]{Remark}
\newtheorem{ntt}[thm]{Notation}

\newcommand{\qed}{\hspace*{\fill}Q.E.D.}  
\newcommand{\cd}{\cdots}
\newcommand{\lr}{\longrightarrow}

\pagenumbering{arabic}


\newcommand{\OA}[1]{{\cal O}_{#1}}
\newcommand{\QA}[1]{{\cal Q}_{#1}}
\newcommand{\SO}[1]{C(S^1) \otimes \OA{#1}}
\newcommand{\sO}[1]{S \otimes \OA{#1}}
\newcommand{\So}[1]{C_0 (S^1 \setminus \{1\}) \otimes \OA{#1}}
\newcommand{\KSO}[2]{K_{#1}(C(S^1) \otimes \OA{#2})}
\newcommand{\KsO}[2]{K_{#1}(S \otimes \OA{#2})}
\newcommand{\KO}[2]{K_{#1} (\OA{#2})}
\newcommand{\Ext}{{\rm Ext}}
\newcommand{\Hom}{{\rm Hom}}
\newcommand{\Tor}{{\rm Tor}}
\newcommand{\ext}[2]{\Ext^1_{{\bf Z}}({#1},{#2})}
\newcommand{\Class}{\mbox{\boldmath ${\cal C}$}}


\newcommand{\CA}{$C^*$-algebra}
\newcommand{\CSalg}{Cuntz-circle algebra}
\newcommand{\DCSalg}{direct limit of Cuntz-circle algebras}
\newcommand{\aab}{approximately absorbing}
\newcommand{\aue}{approximate unitary equivalence}
\newcommand{\ayue}{approximately unitarily equivalent}
\newcommand{\mops}{mutually orthogonal projections}
\newcommand{\hm}{homomorphism}
\newcommand{\pisca}{purely infinite simple \CA}
\newcommand{\aRs}{approximate Rokhlin system}
\newcommand{\andeqn}{\,\,\,\,\,\, {\rm and} \,\,\,\,\,\,}
\newcommand{\QED}{\rule{1.5mm}{3mm}}
\newcommand{\smspace}{\vspace{0.6\baselineskip}}
\newcommand{\bigspace}{\vspace{\baselineskip}}


\newcommand{\arrow}{\rightarrow}
\newcommand{\tdsum}{\widetilde{\oplus}}
\newcommand{\pa}{\|}  
\renewcommand{\ep}{\varepsilon}
\newcommand{\id}{{\rm id}}
\newcommand{\aueeps}[1]{\stackrel{#1}{\sim}}
\newcommand{\aeps}[1]{\stackrel{#1}{\approx}}
\newcommand{\dirlim}{\displaystyle \lim_{\longrightarrow}}

\pagebreak

\begin{abstract}

We prove that if two \hm s from $\OA{\infty}$ to a purely infinite
simple \CA\  have the same class in $KK$-theory, and if either both are
unital or both are nonunital, then they are approximately unitarily
equivalent. It follows that $\OA{\infty}$ is classifiable in the sense
of R\o rdam. In particular, R\o rdam's classification theorem for
direct limits of matrix algebras over
even Cuntz algebras extends to direct limits involving
both matrix algebras over even Cuntz algebras and corners of
$\OA{\infty}$, for which the
$K_0$-group can be an arbitrary countable abelian group with no
even torsion.

\end{abstract}

\vspace{\baselineskip}
\vspace{\baselineskip}

\section{Introduction}

\vspace{\baselineskip}

In \cite{Rr1}, R\o rdam proved that a simple direct limit
$A = \dirlim A_n$ of finite direct sums of matrix algebras over
even Cuntz algebras is classified up to isomorphism by
the pair  $(K_0 (A), [1_A])$ consisting of its
$K_0$-group together with the class of the identity. He furthermore
proved that any pair $(G, g),$ consisting of a countable abelian
odd torsion group $G$ and an element $g \in G,$ can occur
as $(K_0 (A), [1_A])$ for such an algebra $A$.
In this paper, we extend his classification by allowing, as
additional summands in the algebras $A_n,$ matrix algebras
over the infinite Cuntz algebra
$\OA{\infty}$ and arbitrary corners in it. One of the differences
between $\OA{n}$ and $\OA{\infty}$ is that $K_0(\OA{n})$ is
finite while $K_0 (\OA{\infty})$ is infinite cyclic.
This gives a class of
algebras $A$ for which $K_0 (A)$ can be an arbitrary countable
abelian group containing no even torsion, and in which $[1_A]$
can be an arbitrary element of this group.

R\o rdam has defined in \cite{Rr3} a ``classifiable class''
of \pisca s, and has shown that algebras in this class can have
arbitrary countable abelian groups as their $K_0$-groups
(as well as possibly nontrivial $K_1$-groups). Algebras in his
class are determined up to isomorphism by their $K$-theory
together with the class of the identity. The new element in our
work is that we use the natural choice of a \CA\  $A$ satisfying
$K_0 (A) \cong {\bf Z}$ and $K_1 (A) = 0,$ namely $A = \OA{\infty},$
rather than the somewhat arbitrary construction of \cite{Rr3}.
Our results show that our algebras, in particular $\OA{\infty},$
are in fact in R\o rdam's class. It follows that they are isomorphic
to the \CA s with the same $K$-theory constructed in \cite{Rr3}.
Perhaps more importantly, the \CA\  $A$, constructed according to
the recipe in \cite{Rr3} to satisfy $K_0 (A) \cong {\bf Z}$
with generator $[1_A]$ and $K_1 (A) = 0,$ is actually isomorphic
to $\OA{\infty}.$

Our main technical result is that if $D$ is a \pisca , and
if $\varphi, \, \psi : \OA{\infty} \to D$ are two unital
\hm s with the same class in $KK^0 (\OA{\infty}, D),$ then
$\varphi$ is \ayue\  to $\psi.$ Combined with an easy existence
result, this yields the statement that $\OA{\infty}$ is
in R\o rdam's class. (We actually use the simpler Definition 5.1 of
\cite{ER} rather than the definition in \cite{Rr3}.)
The results described above then follow from \cite{Rr3} and a
variation of arguments from \cite{Rr1}.

In the first section, we establish terminology and notation, and
prove several preliminary results, including \aue\  of \hm s
with trivial classes in $KK$-theory. In Section 2, we show that
an arbitrary \hm\  from $\OA{\infty}$ to $D$ is \aab\  in the
sense of \cite{LP} (see Definition 14). The last section contains
the proof of the theorem for \hm s with arbitrary $KK$-classes,
and the consequences discussed above.

This work was done while the first author held a visiting
position at SUNY Buffalo. Some of it was done while both authors were
visiting the Fields Institute. The authors are grateful to both
institutions for their hospitality.

\vspace{\baselineskip}

\section{Preliminaries}

\vspace{\baselineskip}

We begin this section by recalling the definition and standard
properties of the algebras $\OA{\infty}$ and $E_n,$ and establishing
notation for their standard generators. We then define \ayue\  and
\aab\  \hm s. Finally, we prove that if $D$ is a \pisca , then
up to \aue\  there is only one unital \hm\  from $\OA{\infty}$ to
$D$. This is an easy and important special case of our main technical
theorem. We actually work with injective \hm s from
$E_n$ instead, but the two versions of the statement are equivalent.

\vspace{0.6\baselineskip}

{\bf 1.1. Notation.}
Let $\OA{n}$ be the Cuntz algebra, and call its standard generators
$s_1, \dots, s_n.$ Thus the $s_j$ are isometries satisfying
$\sum_{j = 1}^n s_j s_j^* = 1.$

Let $E_n$ be the universal unital \CA\  on generators
$t_1, \dots, t_n$ and relations stating that they are isometries with
orthogonal ranges, but whose range projections do not necessarily
sum to $1.$ Let $J_n$ be the ideal generated by
$1-\sum_{j=1}^nt_jt_j^*.$  It is known (Proposition 3.1 of \cite{Cu1})
that
$J_n \cong \cal K,$ the algebra of compact operators on a
separable infinite dimensional Hilbert space. Let
$\pi_n : E_n \to \OA{n}$ be the quotient map. Thus
$\pi_n (t_j) = s_j,$ and we have a short exact sequence
$$
0 \longrightarrow J_n \longrightarrow E_n
      \stackrel{\pi_n}{\longrightarrow} \OA{n} \longrightarrow 0.
$$
It is known (Proposition 3.9 of \cite{Cu2}) that
$K_0 (E_n) \cong \bf Z$, generated by
$[1]$, and that $K_1 (E_n) = 0$.

The algebra $\OA{\infty}$ can then be viewed as $\dirlim E_n$, where
the map $E_n \to E_{n + 1}$ of the system sends $t_j$ to $t_j$ for
$1 \leq j \leq n.$ Accordingly, we will denote the generators
of $\OA{\infty}$ by $t_1, t_2, \dots ,$ and identify $E_n$ with the
corresponding subalgebra of $\OA{\infty}.$ Recall
(Corollary 3.11 of \cite{Cu2}) that
$K_0 (\OA{\infty}) \cong \bf Z$, generated by
$[1]$, and that $K_1 (\OA{\infty}) = 0$.

\vspace{0.6\baselineskip}

{\bf 1.2 Lemma}. Let $A$ be either $E_n$ or $\OA{\infty}$, and let $D$
be any separable \CA . Then the Kasparov product
$\af \mapsto [1_A] \times \af,$ from $KK^0 (A, D)$ to $K_0 (D),$ is
an isomorphism.

\vspace{0.6\baselineskip}

{\em Proof:} The Universal Coefficient Theorem (\cite{RS}, Theorem 1.18)
shows that this is true for any \CA\  $A$ in the bootstrap category
$\cal N$ of \cite{RS} such that $K_0 (A) \cong \bf Z$ with generator
$[1]$ and $K_1 (A) = 0.$ Thus, we only need to show that our
algebras $A$ are in $\cal N.$
Now $\OA{n}$ is stably isomorphic to a crossed product of an AF
algebra by $\bf Z$ (\cite{Cu1}, 2.1) and so is in $\cal N$.
Since $E_n$ is an extension of $\OA{n}$ by $\cal K$, it too is in
$\cal N$. Therefore $\OA{\infty} \cong \dirlim E_n$ is also in $\cal N$.
\QED

\vspace{0.6\baselineskip}

{\bf 1.3 Lemma}
The standard defining relations for $\OA{n}$ and $E_n$ are exactly
stable in the sense of Loring \cite{Lr4}. That is:

(1) For each $\dt > 0,$ let $\OA{n} (\dt)$ be the universal unital
\CA\  on generators $s_{j, \dt}$ for $1 \leq j \leq n$ and relations
\[
\| s_{j, \dt}^* s_{j, \dt} - 1 \| \leq \dt  \andeqn
     \|\sum_{k = 1}^n s_{k, \dt} s_{k, \dt}^* - 1 \| \leq \dt,
\]
and let $\kp_{\dt} : \OA{n} (\dt) \to \OA{n}$ be the
\hm\  given by $\kp_{\dt} (s_{j, \dt}) = s_j.$ Then for every $\ep > 0$
there is $\dt > 0$ such that there is a
\hm\  $\varphi_{\dt} : \OA{n} \to \OA{n} (\dt)$ satisfying
$\kp_{\dt} \circ \varphi_{\dt} = \id_{\OA{n}}$ and
$\| \varphi_{\dt} (s_j) - s_{j, \dt} \| < \ep$ for all $j.$

(2) For each $\dt > 0,$ let $E_n (\dt)$ be the universal unital
\CA\  on generators $t_{j, \dt}$ for $1 \leq j \leq n$ and relations
\[
\| t_{j, \dt}^* t_{j, \dt} - 1 \| \leq \dt
                      \andeqn
     \| (t_{j, \dt} t_{j, \dt}^*) ( t_{k, \dt} t_{k, \dt}^*) \| \leq \dt
\]
for $j \neq k,$
and let $\kp_{\dt} : E_n (\dt) \to E_n$ be the
\hm\  given by $\kp_{\dt} (t_{j, \dt}) = t_j.$ Then for every $\ep > 0$
there is $\dt > 0$ such that there is a
\hm\  $\varphi_{\dt} : E_n \to E_n (\dt)$ satisfying
$\kp_{\dt} \circ \varphi_{\dt} = \id_{E_n}$ and
$\| \varphi_{\dt} (t_j) - t_{j, \dt} \| < \ep$ for all $j.$

\vspace{0.6\baselineskip}

{\it Proof:}
Part (1) has already been observed to follow from results in the
literature; see the proof of Lemma 2.1 of \cite{LP}.
In any case, it can be proved directly by the same argument as
for (2). We therefore only do (2).  Since the methods are standard, we
will be somewhat sketchy. In particular, we will successively
form the $-1/2$ powers of
elements $x_0, \dots, x_{n - 1}$ of $E_n (\dt)$ (namely the elements
$((1-q_k)t_{k + 1, \dt})^* (1-q_k)t_{k + 1, \dt}$ below).
Each will be a perturbation of $t_{k + 1, \dt}^* t_{k + 1, \dt}$
satisfying
$\|x_k - 1 \|$ small for $\dt$ small enough, but for fixed $\dt$
the norms
$\|x_k - 1 \|$ will grow fairly rapidly with $k.$ The number $\dt$ must
be small enough that all of the $\|x_k - 1 \|$ are small.

Define $w_1 \in E_n (\dt)$ by
$w_1 = t_{1, \dt} (t_{1, \dt}^* t_{1, \dt})^{- 1/2}.$ Then $w_1$ is an
isometry, $\|w_1 - t_{1, \dt} \|$ is small, and $\kp_{\dt} (w_1) = t_1.$

Now suppose we have constructed isometries $w_1, \dots, w_k$
(with $k < n$) such that
$\| w_j - t_{j, \dt} \|$ is small, $\kp_{\dt} (w_j) = t_j$, and
$w_1 w_1^*, \dots, w_k w_k^*$ are \mops . Let
$q_k = w_1 w_1^* + \cdots + w_k w_k^*.$  Since $w_j w_j^*$ is
close to $t_{j, \dt} t_{j, \dt}^*$, it follows from the relations
for $E_n (\dt)$ that $\|q_k t_{k + 1, \dt} t_{k + 1, \dt}^* q_k\|$
is small.
Therefore so is $\|t_{k + 1, \dt}^* q_k t_{k + 1, \dt} \|$, whence
$(1-q_k) t_{k + 1, \dt}$ is close to $t_{k + 1, \dt}.$  Now set
\[
w_{k + 1} = (1-q_k) t_{k + 1, \dt}
              [((1-q_k)t_{k + 1, \dt})^* (1-q_k)t_{k + 1, \dt}]^{-1/2}.
\]
It follows that $w_{k + 1}$ is an isometry which is
close to $t_{k + 1, \dt}$, has range orthogonal to $q_k$, and satisfies
$\kp_{\dt} (w_{k + 1} ) = t_{k + 1}.$

Since only finitely many steps are required, if $\dt$ is small
enough we obtain by induction
isometries $w_1, \dots, w_n \in E_n (\dt)$ with orthogonal
ranges such that $\kp_{\dt} (w_j) = t_j.$ If $\dt$ is sufficiently
small, then we will also have $\|w_j - t_{j, \dt} \| < \ep.$
We then define $\varphi_{\dt} (t_j) = w_j.$ \QED

\vspace{0.6\baselineskip}

The following three definitions are essentially the same as
corresponding definitions in \cite{LP}.

\vspace{0.6\baselineskip}

{\bf 1.4 Definition}. Let $A$ and $B$ be \CA s, let $F$ be a
finite subset of
$A$, and let $\varphi$ and $\psi$ be two \hm s from $A$ to
$B$. We say that $\varphi$ and $\psi$ are {\em \ayue\  to within} $\ep$,
with respect to $F$, if there is a unitary $v \in \tilde{B}$ such that
\[
\| \varphi(f) - v \psi(f) v^* \| < \ep
\]
for all $f \in F$. We abbreviate this as
\[
\varphi \aueeps{\ep} \psi
\]
with respect to $F.$ When the set $F$ is understood, we omit mention
of it.

We further say that $\varphi$ and $\psi$ are {\em \ayue}\  if for every
finite $F \subset A$ and $\ep > 0,$ we have $\varphi \aueeps{\ep} \psi$
with respect to $F.$

\vspace{0.6\baselineskip}

We make two comments on this definition. First, if $A$ has a finite
generating set $G,$ and if $\varphi \aueeps{\ep} \psi$
with respect to $G$ for all $\ep > 0,$ then $\varphi$ is \ayue\  to
$\psi.$ Second, if $F_1 \i F_2 \i \cdots$ are finite subsets of $A$
such that $G = \bigcup_{n = 1}^{\infty} F_n$ generates $A,$ and if
$\varphi \aueeps{\ep} \psi$ with respect to $F_n$ for all $n$ and all
$\ep > 0,$ then again $\varphi$ is \ayue\  to $\psi.$

\vspace{0.6\baselineskip}

{\bf 1.5. Definition} Let $A$ be any unital \CA ,
 and let $D$ be a \pisca .
 Let $\varphi, \psi: A \to D$ be two homomorphisms, and assume that
$\varphi(1) \neq 0$ and
$[\psi(1)] = 0$ in $K_0 (D)$. We define a homomorphism
$\varphi \tdsum \psi : A \to D$, well defined up to unitary equivalence,
by the following construction. Choose a projection $q \in D$ such that
$0 < q < \varphi(1)$ and $[q] = 0$. Since $D$ is purely infinite and
simple, there are partial isometries $v$ and $w$ such that
$  vv^* = \varphi(1) -q$, $v^*v = \varphi(1)$, $ww^* = q$,
 and $w^*w = \psi(1)$.
Now define $(\varphi \tdsum \psi)(a) = v\varphi(a)v^* + w\psi(a)w^*$
for $ a \in A$.

\vspace{0.6\baselineskip}

{\bf 1.6  Definition} Let $D$ be a \pisca , let $A$ be $\OA{n},$
$\OA{\infty},$  $E_n,$  or a finite matrix algebra over one of these,
 and let
 $\varphi: A \to D$ be a homomorphism.
Then $\varphi$ is {\em \aab} if for every $\psi : A \to D$ such that
$[\psi] = 0$ in $KK^0 (A, D),$ the
\hm s $\varphi$ and $\varphi \tdsum \psi$ are \ayue .

\vspace{0.6\baselineskip}

The following proposition is stated in terms of $E_n,$ but it
immediately
implies the same statement about \hm s from $\OA{\infty}.$

\vspace{0.6\baselineskip}

{\bf 1.7 Proposition} Let $D$ be a unital \pisca . Let
$\varphi,\,\psi: E_n\to D$ be two injective
homomorphisms, either both unital or both nonunital.
If $[\varphi]=[\psi]=0$ in $KK^0 (E_n,D),$
then $\varphi$ and $\psi$ are \ayue.

\vspace{0.6\baselineskip}

{\it Proof:} We will  show that, for any $\ep>0$
and any integer $n,$ there is a unitary $w\in D$ such that
$$
\|w^*\varphi(t_j)w-\psi(t_j)\|<\ep
$$
for $j=1,2, \dots,n.$
Since $[\varphi(1)]=[\psi(1)],$ by applying an inner automorphism,
we may assume that
$\varphi(1)=\psi(1).$ Furthermore, replacing $D$ by
$\varphi(1)D\varphi(1),$ we may assume that $\varphi(1)=1.$
For each $j,$ we have
$$
 [\varphi(t_jt_j^*)]=[\psi(t_jt_j^*)]
$$
in $K_0(D).$ Since $\{\varphi(t_jt_j^*)\}$ and
$\{\psi(t_jt_j^*)\}$ are both sequences of mutually orthogonal
nonzero projections in $D,$  there
is a unitary $u\in D$ such that
$$
u^*\varphi(t_jt_j^*)u=\psi(t_jt_j^*)
$$
for $j=1,2, \dots,  n.$ Replacing $\varphi$ by $u^* \varphi (-) u,$
we may assume that
$$
\varphi(t_jt_j^*)=\psi(t_jt_j^*)
$$
for $j=1,2, \dots, n.$ Let $p_j$ be this common value, and set
$q =1-\sum_{j=1}^n p_j.$ Since $[\varphi]=[\psi]=0,$
we have
$$
[1] = \left[\sum_{j = 1}^n p_j \right] = [q] = 0
$$
in $K_0(D).$
Therefore there are  $g_1, g_2, \dots,  g_n \in D$ such that
$$
g_j^* g_j=1 \andeqn \sum_{j=1}^n g_j g_j^* = q.
$$
Define ${\tilde \varphi}: \OA{2n} \to D$ by
$$
{\tilde \varphi}(s_j)=\varphi(t_j) \andeqn
{\tilde \varphi}(s_{n+j}) = g_j
$$
for $j = 1, 2, \dots,  n.$

Let $v_0 = \sum_{j=1}^n \psi(t_j) \varphi(t_j)^*.$
Then $v_0$ is a unitary
in $(1-q)D(1-q).$ Since $D$ is purely infinite and simple,
there is a unitary
$v_1\in qDq$
such that $[v_1]=[v_0^*]$ in $K_1(D).$
Set $w_j=v_1g_j.$ Then $\sum_{j=1}^nw_jg_j^*=v_1.$ Define
${\tilde \psi}: \OA{2n}\to D$ by
$$
{\tilde \psi}(s_j)={\psi(t_j)}  \andeqn
{\tilde \psi}(s_{n+j})=w_j
$$
for $j=1,2, \dots,  n.$
Then
$$
\sum_{j=1}^{2n} {\tilde \psi} (s_j) {\tilde \varphi} (s_j)^* =
               v_0 + v_1 \in U_0 (D).
$$
By Theorem 3.6 of \cite{Rr1}, there exists a unitary $w \in D$ such that
$$
\|w^* {\tilde \varphi} (s_j) w - {\tilde \psi} (s_j) \| < \ep
$$
for $j=1,2, \dots,  2n.$ In particular,
$$
\| w^* \varphi(t_j) w - \psi(t_j)\| < \ep
$$
for $j=1,2, \dots,  n.$
\QED

\section{Approximate absorption}

\vspace{\baselineskip}

In this section, we prove that if $D$ is a \pisca , then any injective
\hm\  from $E_n$ to $D$ is \aab . (This result immediately implies
a corresponding result for \hm s from $\OA{\infty}$ to $D.$) As we will
see in the next section, \aue\  of \hm s with the same class in
$KK$-theory will follow easily.

The technical part of this section is the construction, for even $n$,
of an approximately central projection in $E_n$ whose $K_0$-class is
zero. We first construct a copy of $\OA{n} \oplus E_n$ inside $E_n$
such that the class of $(1, 0)$ in $K_0 (E_n)$ is trivial. Then
we use R\o rdam's work for even Cuntz algebras, and Voiculescu's
Theorem, to move this copy so that the image of $(s_j, t_j)$ is
close to $t_j.$

We mention another approach, not used here since it takes longer to
write. We really only need an approximately central projection in
$\OA{\infty}.$ The inclusion of $E_n$ in $\OA{\infty}$ can be extended
to an inclusion of a suitable Cuntz-Krieger algebra
$\OA{A}$ in $\OA{\infty}$, in such a way that
$K_0 (\OA{A}) \to K_0 (\OA{\infty})$ is an isomorphism. We have
shown that the shift on a Cuntz-Krieger algebra (see Section 4
of \cite{Rr2}) satisfies a version of the approximate Rokhlin
property of \cite{BEK}. This yields approximately central projections,
which with a little extra work can be chosen to be trivial in
$K_0.$

\vspace{0.6\baselineskip}

{\bf 2.1 Lemma}
Let $p$, $q \in E_n \setminus J_n$ be two projections.
If $[p] = [q]$ in $K_0 (E_n),$ then $p$ is Murray-von Neumann
equivalent to $q.$

\vspace{0.6\baselineskip}

{\it Proof:}
Let $\cal P$ be the set of projections in $E_n$ which are
not in $J_n$.  We show that $\cal P$ satisfies the conditions
$(\Pi_1)$--$(\Pi_4)$ before 1.3 of \cite{Cu2}. The conclusion will
then follow from Theorem 1.4 of \cite{Cu2}.

Conditions $(\Pi_1)$ (the sum of two orthogonal projections in
$\cal P$ is again in $\cal P$), $(\Pi_2)$ ($\cal P$ is closed under
Murray-von Neumann equivalence), and $(\Pi_4)$ ($p \in \cal P$ and
$p \leq q$ imply $q \in \cal P$) are all obvious. We thus prove
$(\Pi_3)$. That is, we must show that if
$p$, $q \in E_n \setminus J_n$ are projections, then there is a
projection
$p' \in   E_n \setminus J_n$ such that $p' < q,$ $p'$ is
Murray-von Neumann equivalent to $p,$ and $q-p'\in E_n\setminus J_n.$
Now note that any projection is Murray-von Neumann
equivalent to a subprojection of $t_1 t_1^*,$ since this projection
is Murray-von Neumann equivalent to $1.$ Therefore it suffices to
prove the condition with $p = t_1 t_1^*.$

Since $\OA{n}$ is
purely infinite and simple, we can find a projection
${\bar f} < \pi_n (q)$ and a unitary ${\bar u} \in \OA{n}$
such that ${\bar u} \pi_n (p) {\bar u}^* = {\bar f}.$ Since
$U( \OA{n})$ is connected, there is $u \in U (E_n)$ such that
$\pi_n (u) = {\bar u}.$ Let $f = u p u^*;$ then
$\pi_n (f) < \pi_n (q)$.

We now want to find $e \in (1 - q)J_n (1-q)$ and a unitary $v$ such that
$q + e > v f v^*$. If $1 - q \in J_n,$ we can take $e = 1 - q$
and $v = 1.$ So assume
$1 - q \not\in J_n.$ Then $(1 - q)J_n (1-q)$ is isomorphic to
$\cal K$, and so has an approximate identity $\{e_k\}$ consisting
of projections. Since $f - qf \in J_n,$ we have
\[
(1 - q - e_k) f = (1 - q - e_k) (f - qf) = f - qf - e_k (f - qf)
    \to 0
\]
as $ n \to \infty.$ Therefore $(q + e_k) f \to f$. It follows,
for large enough $k$, that $f$ is unitarily equivalent to a
subprojection $g$ of $q + e_k$, via a unitary that is close to $1.$
Then $\pi_n (q + e_k - g)$ is close to $\pi_n (q) - {\bar f},$ and
in particular is not zero. Thus $q + e_k - g \not\in J_n.$

It now suffices to show that $q + e_k$ is Murray-von Neumann equivalent
to a (not necessarily proper)
subprojection of $q$. Note that $e_k$ is a finite sum of
minimal projections in $(1 - q)J_n (1-q)$. Since this algebra is
isomorphic to $\cal K$, all minimal projections are equivalent, and
it suffices to show that there exists a nonzero projection
$e_0 \in (1 - q)J_n (1-q)$ such that $q + e_0$ is Murray-von Neumann
equivalent to a subprojection of $q$. Since
$J_n \cong \cal K,$ it is equivalent to show that
there is a nonzero projection  $e \in q J_n q$ such that $q$ is
Murray-von Neumann equivalent to a subprojection of $q - e.$

We now claim that
partial isometries in the quotient $\pi_n(q)\OA{n}\pi_n(q)$ lift to
partial isometries in $qE_nq.$ This follows from Corollary 2.12,
the proof of Lemma 2.8, and Remark 2.9 in \cite{Zh}, since
$q J_n q,$ being isomorphic to ${\cal K},$ has
real rank zero and trivial $K_1$-group. (This is also known to
operator theorists. Further see the proof of Lemma 2.6
of \cite{Ell1}.)

It follows from this claim that there is a partial isometry
$s_0 \in q E_n q$ such that $\pi_n (s_0)$ is a
proper isometry in $\pi_n(q) \OA{n} \pi_n(q)$. Then $q - s_0^* s_0$ is a
finite rank projection in $q J_n q \cong \cal K.$ Furthermore,
$q - s_0 s_0^* \not\in J_n$, but $J_n$ is an essential ideal, so
$(q - s_0 s_0^*) J_n (q - s_0 s_0^*)$ contains projections of
arbitrarily large rank. Therefore $s_0$ can be extended to an
isometry $s$ in $q E_n q.$ For the same reason as for $s_0$, there
are projections in $(q - s s^*) J_n (q - s s^*)$ of
arbitrarily large rank. The existence of the required projection $e$
now follows. This completes the proof.
\QED

\vspace{0.6\baselineskip}

{\bf 2.2. Lemma} Let $n$ be an even number. Then for any $\ep>0$
there exist a projection $f\in E_n\setminus J_n$
and partial isometries
\[
v_1^{(1)},v_2^{(1)}, \dots,  v_n^{(1)}\in fE_nf \andeqn
   v_1^{(2)},v_2^{(2)}, \dots,  v_n^{(2)}\in (1 - f)E_n(1 - f)
\]
such that $[f]=0$ in $K_0(E_n)$,
$$
(v_j^{(1)})^* v_j^{(1)}=f, \,\,\,\,\,\,
     \sum_{k=1}^n v_k^{(1)} (v_k^{(1)})^*=f,
$$
$$
(v_j^{(2)})^* v_j^{(2)}= 1-f, \,\,\,\,\,\,
     \sum_{k=1}^n v_k^{(2)} (v_k^{(2)})^* < 1-f,
$$
and
$\| v_j^{(1)} + v_j^{(2)} - t_j \| < \ep$
for $j = 1, \dots, n.$

\vspace{0.6\baselineskip}

{\it Proof:}
Let $e_j = t_j (1 - t_1 t_1^*) t_j^*.$ Then $e_1, e_2, \dots,  e_n$ are
$n$ mutually orthogonal projections in $E_n\setminus J_n$ whose classes
are $0$ in $K_0 (E_n).$ Let
$e=\sum_{j=1}^n e_j.$ By Lemma 2.1,
there are $z_j^{(1)} \in e E_n e$ such that
$$
(z_j^{(1)})^* z_j^{(1)} = e \andeqn z_j^{(1)} (z_j^{(1)})^* = e_j.
$$
Similarly, there is  $w \in E_n$ such that
$$
w^*w=1\andeqn ww^*=1-e.
$$
Define $z_j^{(2)} = w t_j w^*$ and $z_j = z_j^{(1)} + z_j^{(2)}$
for $j=1,2, \dots,  n.$ Then
define $\psi: E_n\to E_n$ by $\psi(t_j) = z_j.$
Note that $[\psi] = 1$ in $KK^0 (E_n, E_n)$ by Lemma 1.2.
Set $z_j=\psi(s_j)$ and
note that $(1-e)-\sum_{j=1}^n z_j^{(2)}(z_j^{(2)})^*$ is a rank one
projection in $(1-e)J_n (1-e)$ (which is isomorphic to ${\cal K}$).
It follows that the elements $\pi_n (z_j) \in \OA{n}$ satisfy
\[
\pi_n (z_j)^* \pi_n (z_j) = 1 \andeqn
    \sum_{j = 1}^n \pi_n (z_j) \pi_n (z_j)^* = 1.
\]

We now have two \hm s from $\OA{n}$ to $\OA{n},$ given by
$s_j \mapsto \pi_n (z_j)$ and $s_j \mapsto \pi_n (t_j) = s_j.$
They both induce the identity map on $K_0 (\OA{n}).$
Since $K_1 (\OA{n}) = 0,$ the Universal Coefficient Theorem
(\cite{RS}, Theorem 1.18) implies that they have the same class
in $KK$-theory.
Let $\dt > 0$ and $\eta > 0$ be small numbers (to be chosen below;
$\dt$ will depend on $\eta$).
By Theorem 3.6 of \cite{Rr1} there is a unitary $v\in \OA{n}$ such that
$$
\|v^*\pi_n(z_j)v-\pi_n(t_j)\|<\dt/2
$$
for $j=1,2, \dots,  n.$
Since $U(\OA{n})$ is connected,
there is a unitary $u\in E_n$ such that $\pi_n(u)=v.$
Then there are $a_j \in J_n$ such that
$$
\|u^* z_j u - (t_j + a_j)\| < \dt.
$$
for $j=1,2, \dots,  n.$

If $\dt$ is sufficiently small, then by Lemma 1.3(2) there are
isometries $t_j'\in E_n$ with orthogonal ranges
such that $\pi_n (t_j') = \pi_n (t_j)$ and
$$
\|t_j + a_j - t_j'\|<\eta
$$
for $j = 1, \dots, n.$ It follows that
\[
\| t_j' - u^* z_j u \| < \eta + \dt.
\]

Let $H$ be an infinite dimensional separable Hilbert space.
Recall that $J_n \cong {\cal K}(H).$
A standard result in representation theory
yields a (unique) representation $\rho : E_n \to B(H)$ which
extends this isomorphism. Clearly $\rho$ is irreducible, and it is
faithful because $J_n$ is an essential ideal in $E_n$. Define a second
representation $\sm : E_n \to B(H)$ by $\sigma(t_j)=t_j'.$
We are going to use Voiculescu's Theorem, as stated in Arveson's paper
\cite{Ar}, to prove that $\rho$ and $\sm$ are \ayue .

Note that
\[
\left\| 1 - \sum_{j = 1}^n t_j' (t_j')^* -
    u^* \left( 1 - \sum_{j = 1}^n z_j z_j^* \right) u \right\|
                  < n (\eta + \dt),
\]
and recall that $1 - \sum_{j = 1}^n z_j z_j^*$ is a rank one
projection in $J_n.$ If $\eta + \dt$ is small enough, it follows
that
\[
\sm \left(1 - \sum_{j = 1}^n t_j t_j^* \right) =
         \rho \left(1 - \sum_{j = 1}^n t_j' (t_j')^* \right)
\]
is a rank one projection in ${\cal K}(H).$ Since it is not zero,
$\sm$ is a
faithful representation of $E_n$. Let $H_0 = \overline{\sm (J_n) H},$
the essential subspace of $\sm |_{J_n}.$ Note that it is a reducing
subspace for $\sm.$ Since $\sm (1 - \sum_{j = 1}^n t_j t_j^*)$
has rank one, we conclude that $(\sm |_{J_n}) |_{H_0}$ is irreducible,
and hence unitarily equivalent to $\rho |_{J_n}.$ Standard results
in representation theory now imply that $\sm ( - ) |_{H_0}$ is
unitarily equivalent to $\rho.$

We have now verified the hypotheses of Theorem 5(iii) of
\cite{Ar}: $\rho$ and $\sm$ have the same kernel (namely $\{0\}$),
their compositions with the quotient map from $B(H)$ to the Calkin
algebra have the same kernel (namely $J_n$), and the essential
parts are unitarily equivalent. Since $\sm (t_j) = \rho (t_j')$,
that theorem yields a unitary $W\in B(H)$ such that
$$
\|W^* \rho (t_j) W - \rho (t_j') \| < \eta \andeqn
        W^* \rho (t_j) W - \rho (t_j') \in {\cal K}(H)
$$
for $j=1, 2, \dots, n.$ Since $t_j - t_j' \in J_n$ and
$W^* \rho (t_j) W - \rho (t_j') \in {\cal K}(H),$
we obtain $\rho (t_j) - W^* \rho (t_j) W \in {\cal K}(H),$
whence $W^* \rho (t_j) W \in  \rho (E_n).$
Furthermore, the $C^*$-subalgebra generated by
$\{W^* \rho (t_j) W : j=1,2, \dots, n \}$ contains
$W^*{\cal K}(H) W ={\cal K}(H),$ and therefore also contains
$\rho (t_j).$ This set thus generates $ \rho (E_n).$
Since $\rho$ is faithful, it follows that there is an automorphism
$\alpha: E_n \to E_n$ such that $\rho (\alpha(t_j)) = W^* \rho (t_j) W$
for $j = 1, 2, \dots, n.$ We can then combine earlier estimates to
obtain
$$
\| \alpha(t_j) - u^* z_j u\| < 2 \eta + \dt
$$
for $j = 1, 2, \dots, n.$
Set
\[
f = \alpha^{-1} (u^* e u) \andeqn
           v_j^{(i)} = \alpha^{-1} (u^* z_j^{(i)} u)
\]
for $i = 1, 2$ and $j = 1, \dots, n.$ Note that $[f] = 0$ in
$K_0 (E_n).$
If $ 2\eta + \dt \leq \ep,$ then these are the desired elements.
\QED

\vspace{0.6\baselineskip}

{\bf 2.3 Proposition} Let $D$ be a \pisca\  and let $\varphi:E_n\to D$
be a monomorphism.
Then $\varphi$ is \aab.

\vspace{0.6\baselineskip}

{\it Proof:} Replacing $D$ by $\varphi (1) D \varphi (1),$ we may assume
that $D$ and $\varphi$ are unital.
Since $q = \varphi (1 - \sum_{j=1}^n t_j t_j^*)$
is nonzero,
there are $n+1$ mutually orthogonal nonzero projections
$$
p_{n+1},p_{n+2}, \dots,  p_{2n},e\in qDq
$$
and isometries
\[
\tilde{t}_{n+1},\tilde{t}_{n+2}, \dots,  \tilde{t}_{2n}\in D
\]
such that
$\tilde{t}_j \tilde{t}_j^* = p_j$ for $j=n+1,n+2, \dots,  2n.$
Now let $A \subset D$ be the $C^*$-subalgebra generated by
$\varphi(t_j)$ for $j=1,2, \dots,  n$ and $\tilde{t}_j$ for
$j=n+1,n+2, \dots,  2n.$ Then $A$ is isomorphic to $E_{2n}.$
By Lemma 2.2, for any $\ep > 0$ there is a projection $f \in A$
and unital \hm s
$\psi_1 : \OA{2n} \to fDf$ and $\psi_2 : E_{2n} \to (1-f) D (1-f)$
such that
$[f]=0$ in $K_0(A)$ (and hence in $K_0(D)$), and
$$
\|\varphi(t_j) - (\psi_1(s_j) + \psi_2 (t_j))\|<\ep/3
$$
for $1 \leq j \leq n$ and
$$
\|\tilde{t}_j - (\psi_1(s_j) + \psi_2 (t_j))\|<\ep/3
$$
for $n + 1 \leq j \leq 2 n.$
Define $\varphi_1, \, \varphi_2 : E_n \to D$ by
$\varphi_1 (t_j) = \psi_1 (s_j)$ and $\varphi_2 (t_j) = \psi_2 (t_j)$
for $j = 1, \dots, n.$
Note that $[\varphi_1]=0$ in $KK^0 (E_n,D)$ by Lemma 1.2.

Now let $\varphi_0 : E_n \to D$ be any \hm\  with
$[\varphi_0] = 0$ in $KK^0 (E_n, D).$ Without loss of generality, we may
assume $\varphi_0 (1) \leq \varphi_1 (1).$ Then
$\varphi_1 \aueeps{\ep / 3} \varphi_1 \tdsum \varphi_0$ by
Proposition 1.7.  Therefore
\[
\varphi \aueeps{\ep/3}  \varphi_1 + \varphi_2 \aueeps{\ep/3}
  (\varphi_1 \tdsum \varphi_0) + \varphi_2 \aueeps{\ep/3}
  \varphi \tdsum \varphi_0,
\]
so $\varphi \aueeps{\ep} \varphi \tdsum \varphi_0$ as desired.
\QED

\section{Classification and direct limits}

\vspace{\baselineskip}

We start this section by proving our main technical theorem, that
\hm s from $\OA{\infty}$ to \pisca s with the same $KK$-classes
are \ayue . This implies that $\OA{\infty}$ is in the
``classifiable class'' $\cal C$ of \cite{ER}
(a slight modification of the class in \cite{Rr3}).
As discussed in the introduction, we then obtain classification
theorems for direct limits involving even Cuntz algebras and
corners of $\OA{\infty},$ for which the $K_0$-groups can have
elements of infinite order. As an interesting corollary, we prove that
if $p \in \OA{\infty}$ is a projection such that $[p] = - [1]$
in $K_0 (\OA{\infty}),$ then $p \OA{\infty} p \cong \OA{\infty}.$

\vspace{0.6\baselineskip}

{\bf 3.1. Definition } Let
$\varphi: \OA{\infty}\to D$ be a homomorphism.
Let $p_j = \varphi (t_j t_j^*)$ for $j = 1, 2, \dots .$
Define ${\bar \varphi}: \OA{\infty}\to (1-p_1-p_2)D(1-p_1-p_2)$
by
$$
{\bar \varphi}(t_j)=\varphi(t_{j+2})(1-p_1-p_2)
$$
for $j=1,2, \dots .$

\vspace{0.6\baselineskip}

{\bf 3.2 Lemma} Let $\varphi$ and $\bar{\varphi}$ be
as in the previous definition, and set
\[
D_0 = [1\oplus(1-p_1-p_2)]M_2(D)[1\oplus (1-p_1-p_2)].
\]
Then the direct sum
\[
\varphi\oplus{\bar \varphi}: \OA{\infty} \to D_0
\]
satisfies $[\varphi\oplus{\tilde \varphi}] = 0$
in $KK^0 ( \OA{\infty}, D_0).$

\vspace{0.6\baselineskip}

{\it Proof:} Clearly
$[(\varphi\oplus{\tilde \varphi}) (1)] = 0$
in $K_0 (D_0).$ The result is now immediate from Lemma 1.2.
\QED

\vspace{0.6\baselineskip}

{\bf 3.3. Theorem} Let $D$ be a \pisca, and
let $\varphi,\,\psi: \OA{\infty}\to D$ be \hm s
such that $[\varphi]=[\psi]$ in $KK^0 (\OA{\infty},D)$ (equivalently,
$[\varphi (1)] = [\psi (1)]$ in $K_0 (D)$).
If $\varphi$
and $\psi$ are both unital or both nonunital, then
$\varphi$ and $\psi$ are \ayue.

\vspace{0.6\baselineskip}

{\it Proof:} We first reduce to the unital case. If both \hm s
are nonunital, then the hypotheses imply that
$\varphi (1)$ is unitarily equivalent to $\psi (1).$
Conjugating $\psi$ by a suitable unitary, we may
thus assume that $\varphi (1) = \psi (1).$ Now replace
$D$ by $\varphi (1) D \varphi (1).$

We now follow the notation of Definition 3.1.
Clearly there is a partial isometry $W\in M_3(O_{\infty})$ such that
$$
W^*W=1\oplus(1-p_1-p_2)\oplus 1\andeqn WW^*=1 \oplus 0 \oplus 0.
$$
Let $G = \{t_1,t_2, \dots\}$ be the standard (infinite) set of
generators of $\OA{\infty},$ and let
$G_{k}=\{t_1,t_2,  \dots, t_{k}\}$ be the set consisting of the first
$k$ of them. Note that $G_k$ generates the canonical copy of $E_k$
in $\OA{\infty}.$  Both $[\varphi \oplus {\bar \varphi}]$ and
$[{\bar \varphi} \oplus \psi]$
are zero in $KK^0 (\OA{\infty}, D),$ so Lemma 1.2 implies that
$[ (\varphi \oplus {\bar \varphi}) |_{E_k} ] =
                      [ ({\bar \varphi} \oplus \psi) |_{E_k} ] = 0$
in $KK^0 (E_k, D).$ For any $\ep > 0,$ Proposition 2.3 now implies
that
$$
\varphi \aueeps{\ep/2}
   \ W (\varphi \oplus {\bar \varphi} \oplus \psi) W^*
 \aueeps {\ep/2} \psi
$$
with respect to $G_k.$ Since $G$ is the increasing union of the $G_k,$
this implies that $\varphi$ and $\psi$ are \ayue.  \QED

\vspace{0.6\baselineskip}

We can now extend R\o rdam's classification theorem
(from Section 7 of \cite{Rr1}) for direct limits
of even Cuntz algebras. For simplicity,
we consider only the case of simple \CA s. We do have to make one
modification in his setup. Every pair $(G, g),$ in which $G$ is a
cyclic group of odd order and $g$ is an element of $G$, occurs as
$(K_0 (M_k (\OA{m})), [1_{M_k (\OA{m})}])$ for some $k$ and some
even $m$. However, the pair $({\bf Z}, 0)$ does not occur as
$(K_0 (M_k (\OA{\infty})), [1_{M_k (\OA{\infty})}])$ for any $k.$
Therefore we will have to allow corners as well as matrix algebras.
Since $\OA{\infty}$ is purely infinite and simple, every finite matrix
algebra is in fact isomorphic to some corner. To simplify the statements
of the results, we will therefore not consider matrix algebras over
$\OA{\infty}.$

\vspace{0.6\baselineskip}

{\bf 3.4. Theorem} Each nonzero corner $p \OA{\infty} p$ of
$\OA{\infty}$ is in the classifiable class
$\cal C$ of Definition 5.1 of \cite{ER}.

\vspace{0.6\baselineskip}

{\it Proof:}
Let $D$ be a \pisca . In the notation of \cite{ER},
$H( p \OA{\infty} p , D)$
is the group of \aue\  classes of nonzero \hm s from
$ p \OA{\infty} p  \otimes {\cal K}$ to
$D \otimes {\cal K}$ and
$KL( p \OA{\infty} p , D)$ is a certain quotient of
$KK^0 ( p \OA{\infty} p , D).$
We have to prove that the \hm\  from
$H( p \OA{\infty} p , D)$ to
$KL( p \OA{\infty} p , D)$ is bijective.
The group
$KL( p \OA{\infty} p , D)$ is defined after Lemma 5.3 in \cite{Rr3},
and in the case at hand is just
$KK^0 ( p \OA{\infty} p , D),$
since the ${\rm Ext}$ terms in the Universal Coefficient Theorem are
zero.
Since
$ p \OA{\infty} p  \otimes {\cal K} \cong \OA{\infty} \otimes {\cal K}$
and
$KL( p \OA{\infty} p , D) \cong KL(\OA{\infty}, D),$
we may assume $p = 1.$

Let $\{e_{ij} : 1 \leq i, j < \infty\}$ be a complete system of
matrix units for $\cal K$. We have to prove that for any
$\eta \in K_0 (D) \cong KK^0 (\OA{\infty}, D),$ there is
up to \aue\  exactly one nonzero
\hm\  $\varphi : \OA{\infty} \otimes {\cal K} \to D \otimes {\cal K}$
such that $[\varphi ( 1 \otimes e_{11})] = \eta$ in $K_0 (D)$.

For existence, choose a nonzero projection $p \in D$ such that
$[p] = \eta.$
Choose a proper isometry
$v_1 \in pDp,$ then choose an isometry  $v_2 \in pDp$ whose range
projection is a proper subprojection of $p - v_1 v_1^*,$
an isometry  $v_3 \in pDp$ whose range
projection is a proper subprojection of $p - v_1 v_1^* - v_2 v_2^*,$
etc., by induction. Define
$\varphi_0 : \OA{\infty} \to D$ by $\varphi_0 (t_j) =  v_j,$ and
take $\varphi = \varphi_0 \otimes \id_{\cal K}.$

For uniqueness, let
$\varphi, \, \psi : \OA{\infty} \otimes {\cal K} \to D \otimes {\cal K}$
be nonzero \hm s with the same class in $KK$-theory.
Identify $M_n \subset {\cal K}$ with
\[
(e_{11} + \cdots + e_{nn}) {\cal K} (e_{11} + \cdots + e_{nn}).
\]
It suffices to prove that for each $n$, the restrictions
of $\varphi$ and $\psi$ to the corner $\OA{\infty} \otimes M_n$ are
\ayue . Now
$\varphi (1 \otimes \sum_{i = 1}^n e_{ii})$
and
$\psi (1 \otimes \sum_{i = 1}^n e_{ii})$
have the same class in $K_0 (D),$ so are unitarily equivalent.
Therefore we may assume they are equal.
Also, $\varphi (1 \otimes e_{11})$ and
$\psi (1 \otimes e_{11})$ have the same class in $K_0 (D),$ so
there is $v_0 \in D$ such that
\[
v_0 v_0^* = \varphi (1 \otimes e_{11})
     \andeqn
v_0^* v_0 = \psi (1 \otimes e_{11}).
\]
Define $v \in U ((D \otimes {\cal K})^+)$ by
\[
v = 1 - \varphi \left(1 \otimes \sum_{i = 1}^n e_{ii} \right) +
\sum_{i = 1}^n \varphi (1 \otimes e_{i1}) v_0 \psi (1 \otimes e_{1i}).
\]
Then
\[
v \psi (1 \otimes e_{ij}) v^* = \varphi (1 \otimes e_{ij})
\]
for $1 \leq i, j \leq n.$
Therefore, without loss of
generality, we may assume that
\[
\psi (1 \otimes e_{ij}) = \varphi (1 \otimes e_{ij})
\]
for $1 \leq i, j \leq n.$
Now it suffices to prove that
$\varphi |_{\OA{\infty} \otimes {\bf C}e_{11}}$
is \ayue\  to
$\psi |_{\OA{\infty} \otimes {\bf C}e_{11}},$
as \hm s from $\OA{\infty}$ to
$\varphi (1 \otimes e_{11}) D \varphi (1 \otimes e_{11}).$ This
follows from Theorem 3.3.
\QED

\vspace{0.6\baselineskip}

{\bf 3.5. Theorem} Let $A = \dirlim A_n$ and $B = \dirlim B_n$
be two simple direct limits, in which each $A_n$ and each $B_n$
is a finite direct sum of matrix algebras over even Cuntz
algebras $\OA{2k}$ and corners in $\OA{\infty}.$

(1) Suppose that $A$ and $B$ are unital, and that there is an
isomorphism
$$
\alpha: (K_0(A),[1_A])\to (K_0(B),[1_B]).
$$
Then there is an isomorphism $\varphi: A\to B$ such that
$\varphi_*=\alpha.$

(2) Suppose that $A$ and $B$ are nonunital, and that there is an
isomorphism
$$
\alpha: K_0(A) \to K_0(B).
$$
Then there is an isomorphism $\varphi: A\to B$ such that
$\varphi_*=\alpha.$

\vspace{0.6\baselineskip}

{\it Proof:}
The previous theorem, combined with Theorem 5.9 of \cite{Rr3},
shows that these direct limits are in the class
$\cal C$ of \cite{ER}. (See the remarks after Definition 5.1 of
\cite{ER}.) The result now follows from Theorem 5.7 of \cite{Rr3}.
\QED

\vspace{0.6\baselineskip}

{\bf 3.6 Corollary.}
If $p, \, q \in \OA{\infty}$ are nonzero projections satisfying
$[p] = \pm [q]$ in $K_0 (\OA{\infty}),$ then
$p \OA{\infty} p \cong q \OA{\infty} q.$ In particular,
\[
(1 - t_1 t_1^* - t_2 t_2^*) \OA{\infty} (1 - t_1 t_1^* - t_2 t_2^*)
        \cong \OA{\infty}.
\]
\QED

\vspace{0.6\baselineskip}

This result is of course easy if $[p] = [q],$ but seems to be new
in the case $[p] = - [q].$

\vspace{0.6\baselineskip}

{\bf 3.7 Lemma} Let $A = \bigoplus_{i = 1}^m A_i$ and
$B = \bigoplus_{i = 1}^n B_i$ be finite direct sums, in which each
$A_i$ and each $B_i$ is a finite matrix algebra over an
even Cuntz algebra $\OA{2k}$ (with $k$ depending on $i$)
or a corner in $\OA{\infty}.$ Let
$\omega : K_0 (A) \to K_0 (B)$ be a \hm\  such that
$\omega ([1_A]) = [1_B].$ Then there is a unital
\hm\  $\varphi : A \to B$
such that $\varphi_* = \omega$ and such that each partial map
$\varphi_{ij} : A_i \to B_j$ is nonzero.

\vspace{0.6\baselineskip}

{\em Proof:}
This is essentially done in the proof of Theorem 2.6 of \cite{Rr1}.
We need only one additional fact, namely that if $p \OA{\infty} p$
is a nonzero corner in $\OA{\infty}$, if $D$ is a \pisca , and if
$\omega : K_0 (\OA{\infty}) \to K_0 (D)$ is a \hm\  such that
$\omega ([p]) = [1_D],$ then there exists a unital
\hm\  $\varphi : p \OA{\infty} p \to D.$ (Recall that
$K_0 (\OA{\infty}) \cong {\bf Z},$ generated by $[1],$ so that
necessarily $\varphi_* = \omega.$)
Choose a projection $q \in D$ such that
$[1_D \oplus q] = \omega ([1_{\OA{\infty}}])$ in $K_0 (D),$
with $q = 0$ if $p = 1$ and $q \neq 0$ otherwise.
Construct a unital
\hm\  $\psi : \OA{\infty} \to (1 \oplus q) M_2 (D) (1 \oplus q),$
as in the existence part of the proof of Theorem 3.4.
In $K_0 (D),$ we then have
$[1_D] = \omega ([p]) = [\psi_* (p)]$ (because
$\omega ([1_{\OA{\infty}}]) = [\psi (1_{\OA{\infty}})]$).
Therefore there
is a unitary $u \in (1 \oplus q) M_2 (D) (1 \oplus q)$ such that
$u \psi (p) u^* = 1 \oplus 0.$ Now take
$\varphi = u \psi(-) u^* |_{p \OA{\infty} p},$ regarded as a \hm\  from
$p \OA{\infty} p$ to $D.$
\QED

\vspace{0.6\baselineskip}

{\bf 3.8. Theorem} Let $G$ be a countable
abelian group with no odd torsion.

(1) Let $g\in G.$ Then there is
a unital simple \CA\ $A,$ a direct limit of finite direct sums of
matrix algebras over even Cuntz algebras $\OA{2k}$ and
corners of $\OA{\infty}$ as in Theorem 3.5,
such that $(K_0(A),[1]) \cong (G,g).$

(2) There is
a simple \CA\ $A$ as in (1), except nonunital,
such that $K_0(A) \cong G.$

\vspace{0.6\baselineskip}

{\it Proof:}
We prove only the unital case. Write
$G = \bigcup_{n = 1}^{\infty} G_n,$ where each
$G_n$ is a finitely generated subgroup of $G,$
and $g \in G_1 \subset G_2 \subset \cdots .$ Each $G_n$ is a finite
direct sum of cyclic subgroups of odd or infinite order.
Therefore there is a finite direct sum
$A_n = \bigoplus_{i = 1}^{r(n)} A_{n, i},$ in which
each $A_{n, i}$ either has the form
$M_{k (n, i)} ( \OA{m(n, i)})$ with $m(n, i)$ even,
or has the form $p_{n, i} \OA{\infty} p_{n, i},$ with
$p_{n, i} \in \OA{\infty}$ a nonzero projection, such that
$K_0 (A_n) \cong G_n.$
With suitable choices of $k(n, i)$ and $p_{n, i},$
we can arrange that this isomorphism sends $[1_{A_n}]$ to $g.$
The previous lemma provides unital \hm s
$\varphi_n : A_n \to A_{n + 1},$ with all partial maps
$A_{n, i} \to A_{n + 1, j}$ nonzero, such that the isomorphisms
$K_0 (A_n) \cong G_n$
and
$K_0 (A_{n + 1}) \cong G_{n + 1}$
identify $(\varphi_n)_*$ with the inclusion of $G_n$ in $G_{n + 1}.$
Now set $A = \dirlim A_n.$ The nontriviality of the partial
embeddings at each stage implies that $A$ is simple. This is the
desired algebra.
\QED

\vspace{0.6\baselineskip}

\end{document}